\documentclass[10pt]{article}
\usepackage{color}
\usepackage[square, comma, sort&compress, numbers]{natbib}
\usepackage{amssymb,amsmath,graphicx,acronym,amsfonts}

\newtheorem{corollary}{Corollary}[section]
\newtheorem{theorem}{Theorem}[section]
\newtheorem{lemma}{Lemma}[section]
\newtheorem{definition}{Definition}[section]
\newtheorem{proposition}{Proposition}[section]
\newtheorem{example}{Example}[section]
\newtheorem{assum}{Assumption}[section]
\newtheorem{algo}{Algorithm}[section]
\newtheorem{Remark}{Remark}[section]

\def\bc{\begin{corl}}
\def\bc{\end{corl}}
\def\ba{\begin{algo}}
\def\ea{\end{algo}}
\def\br{\begin{Remark}}
\def\er{\end{Remark}}
\def\bs{\begin{assum}}
\def\es{\end{assum}}
\def\bt{\begin{theorem}}
\def\et{\end{theorem}\vskip 3pt}
\def\bl{\begin{lemma}}
\def\el{\end{lemma}}
\def\ep{\end{proposition}}
\def\bp{\begin{proposition}}
\def\qed{\hfill{$\Box$}\vskip 5pt}
\def\be{\begin{example}}
\def\ee{\end{example}}
\def\bd{\begin{definition}}
\def\ed{\end{definition}}
\def\bc{\begin{corollary}}
\def\ec{\end{corollary}}
\def\proof{\noindent\it Proof. \hspace{1mm}\rm}
\input epsf
\begin{document}
\title{\bf \Large Finding the maximum eigenvalue of a class of tensors
  with applications in copositivity test and hypergraphs}
\author{Haibin Chen\thanks{School of Management Science, Qufu Normal University, Rizhao, Shandong,
China. Email: chenhaibin508@163.com. This author's work was supported by the National Natural Science Foundation of China (Grant No.  11601261).}, \quad
Yannan Chen\thanks{School of Mathematics and Statistics, Zhengzhou University, Zhengzhou, China. E-mail:
ynchen@zzu.edu.cn. This author's work was supported by the National Natural Science Foundation
of China (Grant No. 11401539) and the Development Foundation for Excellent Youth Scholars of
Zhengzhou University (Grant No. 1421315070).}, \quad
Guoyin Li\thanks{Department of Applied Mathematics, University of New South Wales, Sydney 2052, Australia. E-mail: g.li@unsw.edu.au.
This author's work was partially supported by Australian Research Council.}, \quad
Liqun Qi
\thanks{Department of Applied Mathematics, The Hong Kong Polytechnic University, Hung Hom,
Kowloon, Hong Kong. Email: maqilq@polyu.edu.hk. This author's work
was supported by the Hong Kong Research Grant Council (Grant No.
PolyU 501212, 501913, 15302114 and 15300715).}}
\date{}
\maketitle
\vspace{-0.6cm}
\begin{abstract}
Finding the maximum eigenvalue of a symmetric tensor is an important topic in
tensor computation and numerical multilinear algebra. This paper is devoted to a semi-definite program algorithm for computing the maximum
$H$-eigenvalue of a class of tensors with sign structure called $W$-tensors. The class of $W$-tensors extends the well-studied nonnegative tensors and essentially nonnegative
tensors, and covers some important tensors arising naturally from spectral hypergraph theory.
Our algorithm is based on a new structured sums-of-squares (SOS) decomposition result for a nonnegative homogeneous polynomial induced by a $W$-tensor.
This SOS decomposition enables us to show that computing the maximum $H$-eigenvalue of an
even order symmetric $W$-tensor is equivalent to solving a semi-definite program, and hence can be accomplished in polynomial time.
Numerical examples are given to illustrate that the proposed algorithm can be used to find maximum $H$-eigenvalue of an
even order symmetric $W$-tensor with dimension up to $10,000$. We present two applications for our proposed algorithm: we first
provide a polynomial time algorithm for computing the maximum $H$-eigenvalues of large size Laplacian tensors of hyper-stars and hyper-trees;
second, we show that the proposed SOS algorithm can be used to test the copositivity of a multivariate form
associated with symmetric extended $Z$-tensors, whose order may be even or odd.
Numerical experiments illustrate that our structured semi-definite program algorithm is effective and promising.


%


\end{abstract}

\newpage
\section{Introduction}
Finding the extremum (maximum or minimum) eigenvalue of a tensor is an important topic in tensor computation and numerical multilinear algebra.
Various applications of extremum eigenvalues have been found in the current literature \cite{Cooper12,NQZ2009,DQW15,Qi14}.
For example, in \cite{Qi05}, the sign of the minimum (maximum) eigenvalue plays a crucial role in checking the positive semi-definiteness (negative semi-definiteness) of a symmetric tensor which
has applications in stability analysis of nonlinear autonomous systems involved in automatic control.
In magnetic resonance imaging \cite{Qiy2010,CDH-16}, the principal eigenvalues of an even order symmetric tensor associated with the fiber orientation
distribution of a voxel in cerebral white matter denote volume factions of multiple nerve fibers in this voxel. For a connected even-uniform hypergraph, it has been shown that the maximum H-eigenvalues of the Laplacian tensor and the signless Laplacian tensor are equivalent if and only if the hypergraph is odd-bipartite \cite{Hu2015}.

In view of the importance of eigenvalues of tensors, many researchers have devoted themselves to the study of numerical methods for eigenvalues of high order tensors \cite{CS13,Han2013,Hao2015,HCD2015,HLQS,Kolda14,Ni2015}.
In \cite{Cui2014}, Cui et al. proposed a sophisticated Jacobian semi-definite relaxation method, which computes all of the real eigenvalues of a small symmetric tensor.
Generally speaking, it is an NP-hard problem to compute eigenvalues of a tensor even though the involved tensor is symmetric \cite{Lim2013}.
But for some tensors with special structures, large scale problems can be solved. In \cite{Chyn15}, an inexact curvilinear search
optimization method was established to compute extreme eigenvalues of Hankel tensors, whose dimension may up to one million.

Nonnegative tensors is an important class of structured tensors which arises from the study of image science, statistics, and hypergraph theory.
Ng et al. proposed an iterative method for finding the maximum H-eigenvalue of an irreducible nonnegative tensor \cite{NQZ2009}.
However, the NQZ method is not always convergent for irreducible nonnegative tensors.
Liu et al. \cite{Liu2010} improved the NQZ method so that the refined algorithm is always convergent.
Recently, as a more general class than nonnegative tensors, the essentially nonnegative tensors were studied in \cite{HLQS,Zhang13}.
Hu et al. \cite{HLQS} showed that the maximum $H$-eigenvalue of an even order essentially nonnegative tensor can be found by solving a polynomial optimization problem,
which is equivalently reformulated as a semi-definite programming problem. On the other hand, there are also many important classes of structured tensors
which need not to be nonnegative tensors or essentially nonnegative tensors such as the Laplacian tensor of a hypergraph.
This then raise the following natural question: can we compute the maximum $H$-eigenvalue of a given tensor
with possibly negative values on the off-diagonal elements? This is the main motivation of this paper.

Next, we give a sketch of the copositivity of symmetric tensors. The definition of copositive tensors was introduced in \cite{qlq2013}.
Recently, it has been found that copositive tensors have important applications in vacuum stability of a general scalar potential \cite{KK16}, polynomial optimization \cite{Pena14,song2016} and the tensor complementarity problem \cite{CQW, SQ15, SQ}.
With the help of copositive tensors, Kannike \cite{KK16} studied the vacuum stability of a general scalar potential of a few fields, and it is showed that how to find positivity conditions for more complicated potentials.
Pena et al. \cite{Pena14} proved that recent related results for quadratic problems can be further strengthened and generalized to higher order polynomial optimization problems over the cone of completely positive tensors or copositive tensors.
Che, Qi and Wei \cite{CQW} showed that the tensor complementarity problem defined by a strictly copositive tensor has a nonempty and compact solution set. Song and Qi \cite{SQ15} proved that a real tensor is strictly semi-positive if and only if the
corresponding tensor complementarity problem has a unique solution for any nonnegative vector and a real tensor is semi-positive if and only if the corresponding
tensor complementarity problem has a unique solution for any positive vector. It was shown there that a real symmetric tensor is a (strictly) semi-positive tensor if and only if it is
(strictly) copositive. Song and Qi \cite{SQ} further presented global error bound analysis for the tensor complementarity problem defined by a strictly semi-positive tensor.
Thus, copositive and strictly copositive tensors play an important role in the tensor complementarity problem.

In this article, we propose an efficient semidefinite program algorithm to compute the maximum $H$-eigenvalues of
even order symmetric $W$-tensors, which includes nonnegative tensors and essentially nonnegative tensors as a subclass. This algorithm heavily relies on
an important SOS representation result for a nonnegative polynomial induced from a $W$-tensor.
To proceed, we first give the SOS representation result for nonnegativity of a homogeneous polynomial induced
by an even order symmetric $W$-tensor, which implies that the proposed algorithm maybe much more computationally efficient when
the dimension is large and the explicit expression of the subtensors are available.
Another interesting feature of the $W$-tensor is that it also covers the Laplacian tensors  of hyper-star and hyper-tree \cite{Hu2013, Hu2015}, and hence, the maximum $H$-eigenvalues of even order Laplacian tensors  can be efficiently computed by the proposed algorithm.
Numerical examples show that the proposed algorithm can be used to compute the maximum $H$-eigenvalue of some large size
tensors with dimension up to $10,000$. 
Finally, the proposed algorithm is applied to test the copositivity of a multivariate form associated with symmetric extended $Z$-tensors, where
the order of the tensor can be either odd or even.

This paper is organized as follows. In Section 2, we recall the
definitions for tensors and some basic results about homogeneous polynomials.
In Section 3, we propose an algorithm and show that the maximum $H$-eigenvalue of an even order symmetric $W$-tensor can be computed by the proposed algorithm.
In Section 4, we apply the proposed algorithm  to compute the maximum $H$-eigenvalue of a given large size Laplacian tensor  of
hyper-stars and hyper-trees. Numerical examples are also presented to show the efficiency of our method.
In Section 5, we apply the approach to
test the copositivity of a multivariate form associated with symmetric extended $Z$-tensors, where
the order of the tensor can be either odd or even. In Section 6, we present some final remarks and future work.

Before we move on, we make some comments on notation that will be used in the sequel. Let $\mathbb{R}^n$
be the $n$ dimensional real Euclidean space and $\mathbb{R}^n_+$ be the set of all nonnegative vectors in $\mathbb{R}^n$.
The set consisting of all positive integers
is denoted by $\mathbb{N}$. Let $m, n\in \mathbb{N}$ be two natural numbers. Denote $[n]=\{1,2,\cdots,n\}$. Vectors are denoted by bold
lowercase letters ${\bf x}, {\bf y},\ldots$, matrices are denoted by capital letters $A, B, \ldots$, and tensors are written as calligraphic capitals such as
$\mathcal{A}, \mathcal{T}, \ldots.$ The identity tensor $\mathcal{I}$ with order $m$ and dimension $n$ is given by $\mathcal{I}_{i_1\cdots i_m}=1$ if $i_1=\cdots=i_m$ and $\mathcal{I}_{i_1\cdots i_m}=0$ otherwise.
The $i$-th unit coordinate vector in $\mathbb{R}^n$ is denoted by ${\bf e_i}$, $i\in [n]$.

\setcounter{equation}{0}
\section{Preliminaries}
In this section, we collect some basic definitions and facts that will be used later on. Then,
we introduce the definition of $W$-tensors.

 An $m$-th order $n$-dimensional tensor $\mathcal{A}=(a_{i_1i_2\cdots i_m})$ is a multi-array
of entries $a_{i_1i_2\cdots i_m}$, where $i_j \in [n]$ for $j\in [m]$.
If the entries $a_{i_1i_2\cdots i_m}$ are invariant under any permutation of
their indices, then tensor $\mathcal{A}$ is called a symmetric tensor.
The entries $a_{ii\cdots i}, i\in [n]$, are the diagonal entries of $\mathcal{A}$ and the rest are off-diagonal entries.


We note that an $m$-th order $n$-dimensional symmetric tensor $\mathcal{A}$ uniquely determines an $m$-th degree homogeneous
polynomial $f_{\mathcal{A}}({\bf x})$ on $\mathbb{R}^n$: for all ${\bf x}=(x_1,\cdots,x_n)^T
\in \mathbb{R}^n$,
\begin{equation}\label{e21}
f_{\mathcal{A}}({\bf x})= \mathcal{A}{\bf x}^m=\sum_{i_1,i_2,\cdots, i_m\in [n]}a_{i_1i_2\cdots i_m}x_{i_1}x_{i_2}\cdots x_{i_m}.
\end{equation}
Conversely, an $m$-th degree homogeneous polynomial function $f_{\mathcal{A}}({\bf x})$ on $\mathbb{R}^n$ also uniquely corresponds to a symmetric tensor. Furthermore,
an even order tensor $\mathcal{A}$ is called positive semi-definite (positive definite) if
$$f_{\mathcal{A}}({\bf x}) \geq 0 ~(f_{\mathcal{A}}({\bf x})> 0)~~
{\bf for~ all}~~ {\bf x}\in \mathbb{R}^n \backslash \{\bf 0\}.
$$

Recall that for a polynomial $f$ on $\mathbb{R}^n$, we say $f$ is a sums-of-squares (SOS) polynomial if there exist $r \in \mathbb{N}$ and polynomials $f_i$, $i=1,\ldots,r$ such that
$f=\sum_{i=1}^r f_i^2$. Suppose that $m$ is even.  We denote the set consisting of all SOS polynomials of degree $m$ by $\Sigma^2_m[{\bf x}]$. In (\ref{e21}),
if $f_{\mathcal{A}}({\bf x})$ is an SOS polynomial,
then we say tensor $\mathcal{A}$ has an SOS tensor decomposition \cite{chen15, Hu14}. It is clear that a tensor
with SOS tensor decomposition must be a positive semi-definite tensor, but not vice versa.
For all ${\bf x}\in \mathbb{R}^n$, consider a homogeneous polynomial $f({\bf x})=\sum_{\alpha}f_{\alpha}{\bf x}^{\alpha}$ with degree $m$, where $\alpha=(\alpha_1,\cdots,\alpha_n) \in (\mathbb{N} \cup \{0\})^n$, ${\bf x}^{\alpha}=x_1^{\alpha_1}\cdots x_n^{\alpha_n}$ and $|\alpha|:=\sum_{i=1}^n \alpha_i=m$.
Let $f_{m,i}$ be the coefficient associated with $x_i^{m}$. Let ${\bf e_i}$ be the $i$-th unit vector and let
$$\Omega_f=\{\alpha=(\alpha_1,\cdots,\alpha_n) \in (\mathbb{N} \cup \{0\})^n: f_{\alpha} \neq 0 \mbox{ and } \alpha \neq m \, {\bf e_i}, \ i=1,\cdots,n\}.$$
Then, $f$ can always be written as
$$f({\bf x})=\sum_{i=1}^n f_{m,i} x_i^{m}+\sum_{\alpha \in \Omega_f}f_{\alpha}{\bf x}^{\alpha}.$$

We now recall the definitions of eigenvalues and eigenvectors for a tensor \cite{Lim05, Qi05}.

\bd\label{def21}  Let $\mathbb{C}$ be the complex field. Let $\mathcal{A}=(a_{i_1i_2\cdots i_m})$ be a symmetric tensor with order $m$ dimension $n$. A pair $(\lambda, {\bf x})\in \mathbb{C}\times (\mathbb{C}^n\setminus \{{\bf 0}\})$ is called an
eigenvalue-eigenvector pair of tensor $\mathcal{A}$, if they satisfy
$$
\mathcal{A}{\bf x}^{m-1}=\lambda {\bf x}^{[m-1]},
$$
where $\mathcal{A}{\bf x}^{m-1}$ and ${\bf x}^{[m-1]}$ are all n dimensional column vectors given by
$$\mathcal{A}{\bf x}^{m-1}=\left(\sum_{i_2,\cdots,i_m=1}^n a_{ii_2\cdots i_m}x_{i_2}\cdots x_{i_m} \right)_{1\leq i\leq n}$$
and
$${\bf x}^{[m-1]}=(x_1^{m-1},\ldots,x_n^{m-1})^T \in \mathbb{C}^n.$$
\ed

If the eigenvalue $\lambda$ and the eigenvector ${\bf x}$ are real, then
$\lambda$ is called an $H$-eigenvalue of $\mathcal{A}$ and ${\bf x}$ is its corresponding $H$-eigenvector \cite{Qi05}.
An important fact is that
an even order symmetric tensor is positive semi-definite (definite) if and only if all $H$-eigenvalues of the tensor are nonnegative (positive).
It should be noted that even order symmetric tensors always have $H$-eigenvalues.

The following lemma will play an important role in our later analysis \cite{Qi05}.
\bl\label{lema21} Let $\mathcal{A}$ be a symmetric tensor with order $m$ and dimension $n$, where $m$ is even. Denote the
minimum $H$-eigenvalue and maximum $H$-eigenvalue of $\mathcal{A}$ by
$\lambda_{min}(\mathcal{A})$ and $\lambda_{max}(\mathcal{A})$ respectively. Then, we have
$$
\lambda_{min}(\mathcal{A})=\min_{{\bf x\neq {\bf 0}}}\frac{\mathcal{A}{\bf x}^m}{\|{\bf x}\|_m^m}=\min_{\|{\bf x}\|_m=1}\mathcal{A}{\bf x}^m,
~~~\lambda_{max}(\mathcal{A})=\max_{{\bf x\neq {\bf 0}}}\frac{\mathcal{A}{\bf x}^m}{\|{\bf x}\|_m^m}=\max_{\|{\bf x}\|_m=1}\mathcal{A}{\bf x}^m,
$$
where $\|{\bf x}\|_m=\left(\sum_{i=1}^n|x_i|^m\right)^{\frac{1}{m}}$.
\el

Now, we are ready to define $W$-tensors formally.
For $I\subseteq [n]$, we denote by ${\bf x}_I$ the set of variables $\{x_i : i\in I\}$ and by $\mathbb{R}[{\bf x}_I]$ the polynomial ring in these
variables. For a set $S$ with finitely many members, we use $|S|$ to denote its cardinality.
\bd\label{def22}
Let $\mathcal{A}=(a_{i_1i_2\cdots i_m})$ be a tensor with order $m$ and dimension $n$. We say $\mathcal{A}$ is a $W$-tensor if
there exist $s \in \mathbb{N}$ with $s \le n$ and index sets $\Gamma_l \subseteq [n]$, $l\in [s]$ with  $\bigcup_{l=1}^s \Gamma_l=[n]$ and
$\Gamma_{l_1} \neq \Gamma_{l_2}$ for all $l_1 \neq l_2$ such that
\begin{itemize}
\item[{\bf (i)}] either $s=1$ or $|\left(\bigcup_{l=1}^{p-1}\Gamma_l\right)\bigcap \Gamma_p|\leq 1$ for all $2\leq p\leq s$.

\item[{\bf (ii)}] $\mathcal{A}{\bf x}^m=\sum_{l=1}^s \mathcal{A}_{\Gamma_l}{\bf x}_{\Gamma_l}^m$ for all ${\bf x}\in \mathbb{R}^n$,
where for each $l\in [s]$, $\mathcal{A}_{\Gamma_l}$ is a tensor with order $m$ and dimension $|\Gamma_l|$.

\item[{\bf (iii)}] for each $l\in [s]$, the tensor $\mathcal{A}_{\Gamma_l}$ satisfies either one of the following two conditions:

{\rm (1)} there exists $\{\overline{i_1} , \overline{i_2} , \cdots, \overline{i_m}\} \in \Gamma_l$ with $(\overline{i_1} , \overline{i_2} , \cdots, \overline{i_m})\neq (i,\cdots,i)$, $i\in [n]$ such that the off-diagonal entries of $\mathcal{A}_{\Gamma_l}$ equal zero for all $(i_1,i_2,\cdots,i_m) \in \Gamma_l\backslash \{\pi(\overline{i_1} \, \overline{i_2} \cdots \overline{i_m})\},$
 where $\pi(\overline{i_1} \, \overline{i_2} \cdots \overline{i_m})$ denotes all the permutation of $(\overline{i_1} \, \overline{i_2} \cdots \overline{i_m})$.

{\rm (2)} all the off-diagonal entries of $\mathcal{A}_{\Gamma_l}$ are nonnegative.
\end{itemize}
\ed
Recall that a tensor $\mathcal{A}$ is called an essentially nonnegative tensor if its off-diagonal entries, $\mathcal{A}_{i_1,i_2,\cdots,i_m}$ with $\{i_1,i_2,\ldots,i_m\} \notin \{(i,i,\ldots,i): 1 \le i \le n\}$, are all nonnegative.
It is obvious that essentially nonnegative tensors are $W$-tensors, and the converse is not true in general. Moreover, it is easy to check that
condition ${\bf (i)}$ is automatically satisfied if $\Gamma_1,\cdots,\Gamma_s$ are disjoint i.e., $\Gamma_{l_1}\bigcap \Gamma_{l_2}=\emptyset$ for all $l_1\neq l_2$.
Finally, as we will see later in Section 4, the significance of the $W$-tensors is that it not only extends the essentially nonnegative tensors but also
covers important structured tensors which naturally
arises in the hypergraph theory.

\begin{Remark}\label{remark21} From Definition \ref{def22}, it can be verified that if $\mathcal{A}$ is a $W$-tensor and $\mathcal{D}$ is a diagonal tensor,
then $\mathcal{D}+\mathcal{A}$ is also a $W$-tensor.
\end{Remark}

\setcounter{equation}{0}
\section{Maximum $H$-eigenvalue of a symmetric $W$-tensor}
In this section, we show that the maximum $H$-eigenvalue of an even order symmetric $W$-tensor can be computed by solving a semi-definite programming problem, and so,
can be accomplished in polynomial time.
We first recall a useful lemma, which provides us a simple criterion for determining whether a homogeneous polynomial
with only one mixed term is a sum-of-squares polynomial or not \cite{FK11}.

\bl\label{lema31} Assume $b_1,b_2,\cdots,b_n\geq 0$ and $d\in \mathbb{N}$.
Let $a_1,a_2,\cdots,a_n\in \mathbb{N}$ and $\sum_{i=1}^na_i=2d$. Consider the homogeneous
polynomial $f({\bf x})$ defined by
$$f({\bf x})=b_1x_1^{2d}+\cdots+b_nx_n^{2d}-\mu x_1^{a_1}\cdots x_n^{a_n}.$$
Let $\mu_0=2d\prod _{a_i\neq 0, 1\leq i\leq n}(\frac{b_i}{a_i})^{\frac{a_i}{2d}}.$
Then, the following statements are equivalent:
\begin{itemize}
 \item[{\rm (i)}] $f$ is a nonnegative polynomial i.e., $f({\bf x})\geq 0$ for all ${\bf x}\in \mathbb{R}^n$;
 \item[{\rm (ii)}] $f$ is an SOS polynomial.
\end{itemize}
\el

To proceed, we need the following SOS representation result for nonnegativity of a polynomial induced by a symmetric $W$-tensor.
\bt\label{them31}
Let $\mathcal{A}$ be a symmetric $W$-tensor with even order $m$ and dimension $n$ and let $f({\bf x})=-\mathcal{A}{\bf x}^m$.
Let $\mathcal{A}_{\Gamma_l}$ and  $\Gamma_l, l\in [s]$ be defined as in Definition \ref{def22}.
Suppose that $f({\bf x}) \ge 0$ for all ${\bf x} \in \mathbb{R}^n$.
Then, there exist $h_l \in \mathbb{R}[{\bf x}_{\Gamma_l}]$, $l\in [s]$ and $\rho_i^l \in \mathbb{R}$, $i \in [n]$ and $l \in [s]$, such that each $h_l$ is a sum-of-squares polynomial with
$h_l({\bf x}_{\Gamma_l})=-\mathcal{A}_{\Gamma_l}{\bf x}_{\Gamma_l}^m+ \sum_{i \in {\Gamma_l}} \rho_i^l x_i^m$,  for each $i=1,\ldots,n$,  $$\sum_{l \in \Lambda(i)} \rho_i^l =0, \mbox{ with } \Lambda(i)=\{1 \le l \le s: i \in \Gamma_l\},$$
and  
 \[
 f({\bf x}) = h_1({\bf x}_{\Gamma_1}) + \cdots + h_s({\bf x}_{\Gamma_s}).
 \]
\et
\proof As $\mathcal{A}$ is a symmetric $W$-tensor, there exist $s \in \mathbb{N}$ with $s \le n$ and index sets $\Gamma_l \subseteq [n]$, $l\in [s]$ with  $\bigcup_{l=1}^s \Gamma_l=[n]$ and $\Gamma_{l_1} \neq \Gamma_{l_2}$ for all $l_1 \neq l_2$,
such that conditions {\rm (i)}-{\rm (iii)} hold in Definition \ref{def22}. In particular, condition {\rm (ii)} shows that $f=f_1+\cdots+f_s$ with
 $f_l \in \mathbb{R}[{\bf x}_{\Gamma_l}]$, $l\in [s]$, be homogeneous polynomials with degree $m$ given by
\[
f_l({\bf x}_{\Gamma_l)}= -\mathcal{A}_{\Gamma_l} {\bf x}_{\Gamma_l}^m.
\]
Let us prove the conclusion of this theorem by induction on $s$.

{\bf (1)} We first prove the trivial case, i.e. $s=1$. Then condition {\bf (iii)} of Definition \ref{def22} implies that $f({\bf x})$
is either a polynomial with only one mixed term or
\begin{equation}\label{e31}
f({\bf x})=-\mathcal{A}{\bf x}^m=-\sum_{i=1}^n a_{ii\cdots i}x_i^m-\sum_{\delta_{i_1i_2\cdots i_m} \neq 0}a_{i_1i_2\cdots i_m}x_{i_1}\cdots x_{i_m},
\end{equation}
where $\delta_{i_1i_2 \cdots i_m}$ equals one if $i_1=i_2=\cdots=i_m$, and equals zero otherwise; and  $a_{i_1i_2\cdots i_m}\geq 0$ for all $i_1,i_2,\cdots, i_m\in [n]$ with $\delta_{i_1i_2\cdots i_m} \neq 0$. If the first case holds, we obtain that
$f$ is SOS from Lemma \ref{lema31} since $f({\bf x})\geq 0$ for all ${\bf x}\in \mathbb{R}^n$. If (\ref{e31}) is true, then $\mathcal{A}$ is an even order essentially nonnegative tensor.
It follows from \cite[Proposition 3.1]{HLQS} that $f$ is SOS (see also \cite{Hu14}), and hence, the desired result holds in the case $s=1$.

{\bf (2)}[{\bf Initial Step}] Let $1<s \le n$. We start with $s=2$. Then, it holds that $f=f_1+f_2$ where
$$
f_1({\bf x}_{\Gamma_1})=-\mathcal{A}_{\Gamma_1}{\bf x}_{\Gamma_1}^m \mbox{ and } f_2({\bf x}_{\Gamma_2})=-\mathcal{A}_{\Gamma_2}{\bf x}_{\Gamma_2}^m.
$$
From condition {\bf (i)} of Definition \ref{def22}, we see that there exists $i_0\in [n]$ such that $\Gamma_1\bigcap \Gamma_2\subseteq \{i_0\}$.

If $\Gamma_1\bigcap \Gamma_2=\emptyset$, it can be easily verified that $f_l\geq 0, l\in \{1, 2\}$ since $f({\bf x})\geq 0$.
So, condition {\bf (iii)} of Definition \ref{def22} implies that each $f_l$ is either a homogeneous polynomial with only one mixed term or
a homogeneous polynomial such that $f_l({\bf x}_{\Gamma_1})=-\mathcal{A}_l ({\bf x}_{\Gamma_1})^m$ and $\mathcal{A}_l$ is
an even order essentially nonnegative tensor. This means that $f_1, f_2$ are sum-of-squares polynomials. Thus, $f=f_1+f_2$ is SOS
and the desired result follows with $h_1=f_1, h_2=f_2$.

If $\Gamma_1\bigcap \Gamma_2 =\{i_0\}$, denote $\widehat{\Gamma_l}=\Gamma_l\backslash \{i_0\}, l\in \{1, 2\}$. Without loss of generality, we assume
$\widehat{\Gamma_l}\neq\emptyset, l\in \{1, 2\}$ and we order $\Gamma_l$, $l=1,2$, in such a way that
${\bf x}_{\Gamma_1}=({\bf x}_{\widehat{\Gamma_1}},x_{i_0})$ and ${\bf x}_{\Gamma_2}=(x_{i_0},{\bf x}_{\widehat{\Gamma_2}})$. We first see that
\begin{equation}\label{e32}
\alpha=\inf\{f_1({\bf a},1): {\bf a}\in \mathbb{R}^{|\widehat{\Gamma_1}|}\}>-\infty.
\end{equation}
Otherwise, there exists ${\bf a_k}\in \widehat{\Gamma_1}$ such that $f_1({\bf a_k},1)\rightarrow -\infty$. Let ${\bf 1}_{\widehat{\Gamma_2}}\in \mathbb{R}^ {|\widehat{\Gamma_2}|}$ be the vector such that all its entries equal 1. Note that $({\bf a_k},1,{\bf 1}_{\widehat{\Gamma_2}})\in \mathbb{R}^n$, and so
$$
0\leq f({\bf a_k},1,{\bf 1}_{\widehat{\Gamma_2}})=f_1({\bf a_k},1)+f_2(1,{\bf 1}_{\widehat{\Gamma_2}})\rightarrow -\infty,
$$
which is impossible, and hence, (\ref{e32}) is true. Let $q({\bf x})=\alpha x_{i_0}^m$. Define $h_1=f_1-q$ and $h_2=f_2+q$. We now verify that
$$
h_1=f_1-q\geq 0~~{\bf over}~~\mathbb{R}^{|\Gamma_1|}~ ~{\bf and}~~ h_2=f_2+q\geq 0~~{\bf over}~~\mathbb{R}^{|\Gamma_2|}.
$$
To see this, take any $({\bf a}, 1)\in \mathbb{R}^{|\Gamma_1|}$ with ${\bf a}\in \mathbb{R}^{|\widehat{\Gamma_1}|}$, and so
$$
h_1({\bf a}, 1)=f_1({\bf a}, 1)-\alpha=f_1({\bf a}, 1)-\inf\{f_1({\bf a},1):{\bf a}\in \mathbb{R}^{|\widehat{\Gamma_1}|}\}\geq 0.
$$
Noting that $h_1$ is a homogeneous polynomial with an even order degree $m$, it follows that
$$
h_1({\bf x}_{\Gamma_1})\geq 0~~{\bf for~all}~~{\bf x}_{\Gamma_1}=({\bf a},s)\in \mathbb{R}^{|\widehat{\Gamma_1}|}\times \mathbb{R}~~{\bf with}~~s\neq 0.
$$
Then, by the continuity of $h_1$, it shows that $h_1\geq 0$ over $\mathbb{R}^{|\Gamma_1|}$. Moreover, take any $(1, {\bf b})\in \mathbb{R}^{|\Gamma_2|}$
with ${\bf b}\in \mathbb{R}^{|\widehat{\Gamma_2}|}$. Fix an arbitrary $\epsilon>0$. Take ${\bf z}_{\epsilon}\in \mathbb{R}^{|\widehat{\Gamma_1}|}$
be such that $f_1({\bf z}_{\epsilon},1)\leq \inf\{f_1({\bf x}_{|\widehat{\Gamma_1}|},1):{\bf x}_{|\widehat{\Gamma_1}|}\in \mathbb{R}^{|\widehat{\Gamma_1}|}\}+\epsilon$. Then,
$$
\begin{aligned}
h_2(1, {\bf b})=&f_2(1, {\bf b})+\inf\{f_1({\bf x}_{|\widehat{\Gamma_1}|},1):{\bf x}_{|\widehat{\Gamma_1}|}\in \mathbb{R}^{|\widehat{\Gamma_1}|}\} \\
\geq & f_2(1, {\bf b})+f_1({\bf z}_{\epsilon},1)-\epsilon\\
=&f({\bf z}_{\epsilon},1, {\bf b})-\epsilon \geq -\epsilon,
\end{aligned}
$$
where the last inequality follows by the fact that $f({\bf x})\geq 0$ for all ${\bf x}\in \mathbb{R}^n$. Letting $\epsilon \rightarrow 0$, we have
$h_2(1, {\bf b})\geq 0$ for all ${\bf b}\in \mathbb{R}^{|\widehat{\Gamma_2}|}$. Similarly, using the fact that $h_2$ is a continuous homogeneous polynomial with
even degree $m$, we see that $h_2\geq 0$ over $\mathbb{R}^{|\Gamma_2|}$.

To finish the proof of the initial step i.e. $s=2$, it remains to observe from condition {\bf (iii)} of Definition \ref{def22} that
each $h_l\in \mathbb{R}[{\bf x}_{\Gamma_l}]$ is either a polynomial with only one mixed term or a homogeneous polynomial such that $h_l({\bf x}_{\Gamma_l})=-\mathcal{H}_l {\bf {\bf x}_{\Gamma_l}}^m$ and
$\mathcal{H}_l$ is an essentially nonnegative tensor.
Similar to the analysis of {\bf (1)}, we know that $h_l, l\in \{1, 2\}$ are SOS. Therefore, the conclusion holds with $s=2$.

{\bf (3)}[{\bf Induction Step}] Suppose that the conclusion is true for $s=p-1$. We now examine the case with $s=p$.
In this case, we have
$$
f=f_1+f_2+\cdots+f_p=\hat{f}+f_p,
$$
where $\hat{f}=f_1+f_2+\cdots+f_{p-1}$. Then, $\hat{f}$ is a homogeneous polynomial with $\hat{f}\in \mathbb{R}[{\bf x}_{\bigcup_{l=1}^{p-1}\Gamma_l}]$.
By the definition of $W$-tensor, there exists $i_p\in [n]$ such that
$$
\left(\bigcup_{l=1}^{p-1}\Gamma_l\right)\bigcap\Gamma_p\subseteq \{i_p\}.
$$
Similarly to the proof in initial step {\bf (2)}, there exist $\hat{h}\in \mathbb{R}[{\bf x}_{\bigcup_{l=1}^{p-1}\Gamma_l}]$, $h_p\in \mathbb{R}[{\bf x}_{\Gamma_p}]$ and finite number $\rho \in \mathbb{R}$ such that
$$
\hat{h}=\hat{f}-\rho x_{i_p}^m\geq 0~~{\bf over}~~\mathbb{R}^{|\bigcup_{l=1}^{p-1}\Gamma_l|},~  h_p=f_p+\rho x_{i_p}^m\geq0~~{\bf over}~~\mathbb{R}^{|\Gamma_p|}.
$$
By the induction hypothesis that the conclusion holds when $s=p-1$, we know that $\hat{h}=h_1+h_2+\cdots+h_{p-1}$, where
$h_l \in \mathbb{R}[{\bf x}_{\Gamma_l}]$, $l\in [p-1]$, such that each $h_l$ is a sum-of-squares polynomial with
$h_l({\bf x}_{\Gamma_l})=\mathcal{A}_{\Gamma_l}{\bf x}_{\Gamma_l}^m+ \sum_{i \in {\Gamma_l}} \rho_i^l x_i^m$ for some $\rho_i^l \in \mathbb{R}$ satisfying $\sum_{l=1}^{p-1}\sum_{i \in \Gamma_l}\rho_i^l=0$
and  
 \[
 \hat{f}({\bf x}) = h_1({\bf x}_{\Gamma_1}) + \cdots + h_{p-1}({\bf x}_{\Gamma_{p-1}}).
 \]
 On the other hand, $h_p=f_p+\rho x_{i_p}^m$ is either a homogeneous polynomial with one mixed term or a homogeneous polynomial
 a homogeneous polynomial such that $h_p({\bf x})=-\mathcal{H}_p {\bf {\bf x}}^m$ and
$\mathcal{H}_p$ is an essentially nonnegative tensor.
 Thus $h_p$ is SOS and the desired results hold. \qed

\bt\label{them32} Let $\mathcal{A}$ be a symmetric $W$-tensor with even order $m$ and dimension $n$.
Then, it holds that
\begin{equation}\label{e33}
\lambda_{max}(\mathcal{A})= \min_{t\in \mathbb{R},\mu\in \mathbb{R}} \{t~|~t-\mathcal{A}{\bf x}^m+\mu \left(\|{\bf x}\|_m^m-1\right)\in \Sigma_m^2[{\bf x}]\}.
\end{equation}
Moreover, let $\mathcal{A}_{\Gamma_l}, l\in [s]$ be subtensors of $\mathcal{A}$ as defined in Definition \ref{def22}, then we also have
\begin{equation}\label{e33+}
\lambda_{max}(\mathcal{A})=\min_{t,\rho_i^l\in \mathbb{R}}\{t : -\mathcal{A}_{\Gamma_l}{\bf x}_{\Gamma_l}^m+\sum_{i\in \Gamma_l}\rho_i^lx_i^m\in \Sigma_m^2[{\bf x}_{\Gamma_l}], l\in [s], \sum_{l \in \Lambda(i)} \rho_i^l \le t, \,  i \in [n]\},
\end{equation}
where, for each $i=1,\ldots,n,$ $\Lambda(i)=\{1 \le l \le s: i \in \Gamma_l\}$.
\et
\proof Let $t^*=\lambda_{max}(\mathcal{A})$. By Lemma \ref{lema21}, we have that
$$
f({\bf x})=t^*\sum_{i=1}^mx_i^m-\mathcal{A}{\bf x}^m\geq 0,~~{\bf for ~all}~~{\bf x}\in \mathbb{R}^n.
$$
Let $\mathcal{B}=\mathcal{A}-t^*\mathcal{I}$. It can be easily verified that $f({\bf x})=-\mathcal{B}{\bf x}^m\geq 0$ for all ${\bf x}\in \mathbb{R}^n$.
By Definition \ref{def22} and Remark 2.1, we know that $\mathcal{B}$ is a symmetric $W$-tensor since $\mathcal{A}$ is $W$-tensor. Let $\mathcal{B}_{\Gamma_l}, l\in [s]$
denote the subtensors of $\mathcal{B}$. So, it is easy to verify that $\mathcal{B}_{\Gamma_l}$ only differs $\mathcal{A}_{\Gamma_l}$ by a
diagonal tensor on $\mathbb{R}^{|\Gamma_l|}$. By Theorem \ref{them31}, there exist sum-of-squares polynomials $h_l\in \mathbb{R}[{\bf x}_{\Gamma_l}], l\in [s]$
such that
\begin{equation}\label{e34}
f=h_1+h_2+\cdots+h_s
\end{equation}
and
\begin{equation}\label{e35}
h_l({\bf x}_{\Gamma_l})=-\mathcal{A}_{\Gamma_l}{\bf x}_{\Gamma_l}^m+\sum_{i\in \Gamma_l}\bar{\rho}_i^lx_i^m
\end{equation}
for some $\bar{\rho}_i^l\in \mathbb{R}$. Thus, $f({\bf x})=-\mathcal{B}{\bf x}^m=-\mathcal{A}{\bf x}^m+t^* \|x\|_m^m$ is an SOS polynomial, and so, $t=\mu=t^*$ is feasible for the problem
$$
\min_{t\in \mathbb{R},\mu\in \mathbb{R}} \{t~|~t-\mathcal{A}{\bf x}^m+\mu \left(\|{\bf x}\|_m^m-1\right)\in \Sigma_m^2[{\bf x}]\}.
$$
So, it follows that
$$
\min_{t\in \mathbb{R},\mu\in \mathbb{R}} \{t~|~t-\mathcal{A}{\bf x}^m+\mu \left(\|{\bf x}\|_m^m-1\right)\in \Sigma_m^2[{\bf x}]\}\leq t^*=\lambda_{max}(\mathcal{A}).
$$

On the other hand, for all ${\bf x}\in \mathbb{R}^n$, take any $(t,\mu)$ with $t-\mathcal{A}{\bf x}^m+\mu (\|{\bf x}\|_m^m-1)\in \Sigma_m^2[{\bf x}].$
Then, it holds that
$$t-\mathcal{A}{\bf x}^m+\mu (\|{\bf x}\|_m^m-1)\geq 0,~~\forall~{\bf x}\in \mathbb{R}^n,$$
which implies that $t\geq \mu$ and $\mu \|{\bf x}\|_m^m\geq \mathcal{A}{\bf x}^m$ for all ${\bf x} \in \mathbb{R}^n$.
Thus, we know that $t^*\geq \mathcal{A}{\bf x}^m$ for all ${\bf x}\in \mathbb{R}^n$ with $\|{\bf x}\|_m^m=1$. Therefore, (\ref{e33}) follows.

To see the second part, by (\ref{e34}) and (\ref{e35}), we see that
$$
\sum_{l \in \Lambda(i)} \bar \rho_i^l =t^*.
$$
Combining this with the fact that each $h_l$ is a sum-of-squares polynomial, gives us that $t=t^*$ and $\rho_i^l=\bar{\rho}_i^l$ is feasible for
$$
\min_{t,\rho_i^l\in \mathbb{R}}\{t : -\mathcal{A}_{\Gamma_l}{\bf x}_{\Gamma_l}^m+\sum_{i\in \Gamma_l}\rho_i^lx_i^m\in \Sigma_m^2[{\bf x}_{\Gamma_l}], l\in [s], \sum_{l \in \Lambda(i)} \rho_i^l \le t, \,  i \in [n]\},
$$
which implies that
$$
\min_{t,\rho_i^l\in \mathbb{R}}\{t : -\mathcal{A}_{\Gamma_l}{\bf x}_{\Gamma_l}^m+\sum_{i\in \Gamma_l}\rho_i^lx_i^m\in \Sigma_m^2[{\bf x}_{\Gamma_l}], l\in [s], \sum_{l \in \Lambda(i)} \rho_i^l \le t, \,  i \in [n]\}\leq t^*.
$$
Moreover, by a direct computation, the reverse inequality always hold such that
$$
\min_{t,\rho_i^l\in \mathbb{R}}\{t : -\mathcal{A}_{\Gamma_l}{\bf x}_{\Gamma_l}^m+\sum_{i\in \Gamma_l}\rho_i^lx_i^m\in \Sigma_m^2[{\bf x}_{\Gamma_l}], l\in [s], \sum_{l \in \Lambda(i)} \rho_i^l \le t, \,  i \in [n]\}\geq t^*.
$$
Therefore, the desired results hold. \qed
\begin{Remark}\label{remark31}{\bf (Polynomial time solvability of maximum $H$-eigenvalue of the symmetric $W$-tensor)}
As explained in \cite{Hu14}, checking a sum-of-squares polynomial can be equivalently rewritten as a semi-definite
programming problem. Then, Theorem \ref{them32} shows that if $\mathcal{A}$ is a symmetric $W$-tensor with even order, its maximum $H$-eigenvalue can be found
by solving a semi-definite problem, and so can be validated in polynomial time.
\end{Remark}


According to the equality (\ref{e33}), we obtain a basic algorithms for finding the maximum eigenvalue of a $W$-tensor.

\ba\label{algo31} Given $m,n\in \mathbb{N}$ and let $m$ be even.

{\bf 1} Input an $m$-th order $n$-dimensional $W$-tensor $\mathcal{A}=(a_{i_1i_2\cdots i_m})$.

{\bf 2} Let $f_{\mathcal{A}}({\bf x})=\mathcal{A}{\bf x}^m$ and $g({\bf x})=\|x\|_m^m-1=\sum_{i\in [n]}x_i^m-1.$
Solve the following optimization problem
$$
\min_{t\in \mathbb{R},\mu\in \mathbb{R}} \{t~|~t-f_{\mathcal{A}}({\bf x})+\mu g({\bf x})\in \Sigma_m^2[{\bf x}]\}.
$$

{\bf 3} Output $t$.
\ea

When the subtensors $\mathcal{A}_l$ of a $W$-tensor $\mathcal{A}$ can be explicitly exploited, we can have the following refined algorithm
for computing the maximum eigenvalue of a $W$-tensor.

\ba\label{algo32} Given $m,n\in \mathbb{N}$ and let $m$ be even.

{\bf 1} Input an $m$-th order $n$-dimensional $W$-tensor $\mathcal{A}=(a_{i_1i_2\cdots i_m})$ with its subtensors $\mathcal{A}_{\Gamma_l}$, $l \in [s]$ as defined in Definition \ref{def22}.

{\bf 2}
Solve the following optimization problem
$$
\min_{t,\rho_i^l\in \mathbb{R}}\{t : -\mathcal{A}_{\Gamma_l}{\bf x}_{\Gamma_l}^m+\sum_{i\in \Gamma_l}\rho_i^lx_i^m\in \Sigma_m^2[{\bf x}_{\Gamma_l}], l\in [s],
\sum_{l \in \Lambda(i)} \rho_i^l \le t, \,  i \in [n]\},
$$
where $\Lambda(i)=\{1 \le l \le s: i \in \Gamma_l\}$.

{\bf 3} Output $t$.
\ea

\begin{Remark}{\bf (Comparison with Algorithms \ref{algo31} and \ref{algo32})} It is known that checking
a polynomial with even degree $d$ and dimension $n$ is sums-of-squares or not can be equivalently reformulated as a semi-definite programming problem
which can be done via the Matlab Toolbox {\rm YALMIP} \cite{Yalmip1,Yalmip2}.
Thus, Algorithms 3.1 and 3.2 can both be used to compute the maximum $H$-eigenvalue of a $W$-tensor via {\rm YALMIP} and the commonly used SDP solver such as SeDuMi \cite{Sturm99}. Algorithm 3.2 requires the explicit
expression of the subtensors of a $W$-tensor while Algorithm 3.1 does not need this information. On the other hand, in the case where $n$ is large and the explicit
expression of the subtensors are available, Algorithm 3.2 can be much more computationally efficient. Indeed, note that
checking whether a polynomial with even degree $d$ and dimension $n$ is sums-of-squares or not leads to a semi-definite programming problem whose size (the maximum of the number of the variables and
the number of involved linear constraints) is ${n+d \choose d}$.
So, for an $m$-th order $n$-dimensional tensor, Algorithm 3.1 amounts of solving
a semi-definite programming problem with size ${n+m \choose m}$; while Algorithm 3.2
leads to a semi-definite programming problem with  size $s{k+m \choose m}$ where
$k=\max\{|\Gamma_l|: 1 \le l \le s\}$, which is much smaller than ${n+m \choose m}$ if
$n$ is large and $k$ is small. For example, if $s=n/4$, $k=4$ and $m=4$, then ${n+m \choose m}$ is of the order $n^4$; while $s{k+\frac{m}{2} \choose \frac{m}{2}}$
is of the order $n$.
\end{Remark}

Next, we present an example to illustrate that Algorithm \ref{algo32} can be used to compute the maximum $H$-eigenvalue for large size W-tensors (dimension up to $10,000$).

\begin{example}\label{exam31}
Let $n=4k$ with $k \in \mathbb{N}$. Consider the symmetric tensor $\mathcal{A}$ with order $4$ and dimension $n$ where
\[
\mathcal{A}_{1111}=\mathcal{A}_{2222}=\cdots=\mathcal{A}_{nnnn}=n,
\]
\[
\mathcal{A}_{i_1i_2i_3i_4}=-\frac{1}{6}, \mbox{ for all } (i_1,i_2,,i_3,i_4)=\pi(4l-3,4l-2,4l-1,4l), l=1,\cdots,\frac{n}{4},
\]
and $\mathcal{A}_{i_1i_2i_3i_4}=0$ otherwise. Here $\pi(i_1,\cdots,i_4)$ denotes all the possible permutation of $(i_1,\cdots,i_4)$. Clearly $\mathcal{A}$ is not a nonnegative tensor (or an
essentially nonnegative tensor).
The tensor $\mathcal{A}$ corresponds to a unique homogeneous polynomial
$$f_{\mathcal{A}}({\bf x})=\mathcal{A}{\bf x}^m= n(x_1^{4}+\cdots+x_n^4) - 4\sum_{l=1}^{n/4}x_{4l-3}\, x_{4l-2}\, x_{4l-1}\, x_{4l}.$$
Let $\Gamma_l=(4l-3,4l-2,4l-1,4l)$, $l=1,\cdots,\frac{n}{4}$. Then, by Definition \ref{def22}, $\mathcal{A}$ is a $W$-tensor with subtensors $\mathcal{A}_{\Gamma_l}$, $l=1,\ldots,\frac{n}{4}$, defined by $$\mathcal{A}_{\Gamma_l}{\bf x}_{\Gamma_l}^4=n(x_{4l-3}^{4}+x_{4l-2}^4+x_{4l-1}^4+x_{4l}^4) - 4 x_{4l-3}\, x_{4l-2}\, x_{4l-1}\, x_{4l}. $$
Moreover, using geometric mean inequality, we can directly verify that the true maximum $H$-eigenvalue of $\mathcal{A}$ is $\lambda_{\max}^H(\mathcal{A}) = n+1$.

\begin{table}[!htb]
\begin{center}
\caption{Test results for Example \ref{exam31}}
\begin{tabular}{|cc|c|ccc|}\hline
$m$ & $n$ & True $\lambda_{\max}^H(\mathcal{A})$ & Est. $\lambda_{\max}^H(\mathcal{A})$ & YALMIP & SeDuMi \\ \hline
 4  & ~~500 & ~~501 & ~~501.0000 & ~~10.2 & ~~4.3 \\
 4  & ~1000 & ~1001 & ~1001.0000 & ~~29.9 & ~12.6 \\
 4  & ~2000 & ~2001 & ~2001.0000 & ~112.5 & ~26.3 \\
 4  & ~5000 & ~5001 & ~5001.0000 & ~650.8 & ~66.4 \\
 4  & 10000 & 10001 & 10001.0000 & 2164.8 & 136.7 \\   \hline
\end{tabular}
\end{center}
\end{table}

We compute the maximum $H$-eigenvalue of $\mathcal{A}$ using Algorithm 3.2 for
the cases of $n =500$ to $10,000$, where True $\lambda_{\max}^H(\mathcal{A})$ and Est. $\lambda_{\max}^H(\mathcal{A})$
denote the true maximum $H$-eigenvalue and estimated maximum $H$-eigenvalue respectively.
The results are summarized in Table 1.
Obviously, Algorithm 3.2 finds the maximum $H$-eigenvalue of $\mathcal{A}$ exactly.
The CPU-time (measured in seconds) for converting the sums-of-squares problem to SDP via the Matlab toolbox YALMIP
and solving SDP via the commonly used SDP software SeDuMi \cite{Sturm99} is reported in columns YALMIP and SeDuMi, respectively.
When the dimension of the tensor increases, the CPU-time for YALMIP and SeDuMi grows steadily.

\end{example}

\setcounter{equation}{0}
\section{Applications in spectra of hypergraphs}

Throughout this part, unless stated otherwise, a hypergraph means an undirected simple $k$-uniform hypergraph
$G=(V, E)$, where $E\subseteq 2^V$. The elements of $V=V(G)$, which is labeled as $[n]=\{1,2,\cdots,n\}$, are
referred to as vertices and the elements of $E=E(G)$ are called edges. Recall that a simple hypergraph is a hypergraph
where none of its edges is contained within another. We say a hypergraph is $m$-uniform if for every edge $e\in E$, it
holds that $|e|=m$. For a subset $S\subseteq [n]$, we denote by $E_S$ the set of edges $\{e\in E~ | ~S\cap e\neq \emptyset\}$.
For a vertex $i\in V$, we simplify $E_{\{i\}}$ as $E_i$. The cardinality of the set $E_i$ is defined as the degree of the vertex $i$,
which is denoted by $d_i$.

Now, we first introduce several basic definitions that will be studied.
The following definition for Laplacian tensor and signless Laplacian tensor were proposed by Qi in \cite{Qi14}. For other related definitions see \cite{HQ,Ligy2013}.
\bd\label{def41}{\bf (Laplacian and signless Laplacian of hypergraphs)} Let $G=(V,E)$ be an $m$-uniform hypergraph where $V=\{1,2,\cdots,n\}$. The
adjacency tensor of $G$ is defined as the $m$th order $n$ dimensional tensor $\mathcal{A}$ with
$$
a_{i_1i_2\cdots i_m}=\left\{
\begin{array}{cll}
\frac{1}{(m-1)!} & \{i_1,i_2,\cdots,i_m\}\in E, \\
0 & {\bf otherwise}.
\end{array}
\right.
$$
Let $\mathcal{D}$ be an $m$th order $n$ dimensional diagonal tensor with its diagonal elements equal to the degree of vertex $i$, for all $i\in [n]$.
Then $\mathcal{L}=\mathcal{D}-\mathcal{A}$ is the Laplacian tensor of hypergraph $G$, and $\mathcal{Q}=\mathcal{D}+\mathcal{A}$ is the signless Laplacian tensor of hypergraph $G$.
\ed

It can be easily verified that the signless Laplacian tensor of a hypergraph is a nonnegative tensor, and so, is, in particular, a $W$-tensor. On the other hand, the off-diagonal elements of the Laplacian tensor of a hypergraph can be negative. Below, we show that the Laplacian tensors of two important types of hypergraphs are indeed $W$-tensors, and hence, their maximum $H$-eigenvalue can be found in polynomial time via Algorithm 3.2.


\subsection{Laplacian tensors of hyper-stars}
To move on, we first recall the concept of a hyper-star \cite{Hu2015}.
\bd\label{def43}Let $V=\{1,2,\cdots,n\}$ and $E$ is a set of subsets of $V$. Let $G=(V,E)$ be a $m$-uniform
hypergraph. If there is a disjoint partition of the vertex set $V$ as $V=V_0\cup V_1\cup \cdots \cup V_s$
such that $|V_0|=1$ and $|V_1|=|V_2|=\cdots=|V_s|=m-1$, and $E=\{V_0\cup V_i~|~i\in [s]\}$, then $G$ is called a {\bf hyper-star}.
\ed
It is an immediate fact that, with a possible renumbering of the vertices, all the
hyper-stars with the same size are identical.

\bt\label{them41} Let $G=(V, E)$ be a hyper-star. Then its Laplacian tensor is a symmetric $W$-tensor.
\et
\proof Let $\mathcal{A}$ and  $\mathcal{L}$ be the adjacent tensor and Laplacian tensor of $G$ respectively. Then
$\mathcal{L}=\mathcal{D}-\mathcal{A}$, where $\mathcal{D}$ is the diagonal tensor with its diagonal entries $d_i$ i.e., the degree of the vertex $i\in [n]$.
Assume $V=[n]$ and $|E|=s$. Let $V_0=\{i^0\}$ and $V_l=\{i^l_1,i^l_2,\cdots,i^l_{m-1}\}, l\in [s]$. Define $\Gamma_l=\{i^0,i^l_1,i^l_2,\cdots,i^l_{m-1}\}$, $l\in [s]$. It holds that
$\left( \bigcup_{l=1}^{p-1} \Gamma_l\right)\bigcap\Gamma_p=\{i^0\}$, for all $2\leq p\leq s$.
By Definition \ref{def41} and Definition \ref{def43}, we know that
$$
\mathcal{L}_{i_1i_2\cdots i_m}=\left\{
\begin{array}{cll}
d_i & {\bf if}~~i_1=i_2=\cdots=i_m=i \\
-\frac{1}{(m-1)!} & {\bf if}~~\{i_1,i_2,\cdots,i_m\}=V_0\cup V_l~~{\bf for~some}~~l\in [s] \\
0 & {\bf otherwise}
\end{array}
\right.
$$
and
$$
\mathcal{L}{\bf x}^m=\sum_{i=1}^md_ix_i^m-\sum_{\{i_1,\cdots,i_m\}\in E}\frac{1}{(m-1)!}x_{i_1}x_{i_2}\cdots x_{i_m}.
$$
For any $l\in [s]$, define $m$th order $|\Gamma_l|$ dimensional tensors $\mathcal{L}_{\Gamma_l}$ such that
$$
(\mathcal{L}_{\Gamma_l})_{i_1i_2\cdots i_m}=\left\{
\begin{array}{cll}
\frac{d_{i^0}}{s} & {\bf if}~~i_1=i_2=\cdots=i_m=i^0 \\
d_{i^l_j} & {\bf if}~~i_1=i_2=\cdots=i_m=i^l_j,~j=1,2,\cdots,m-1 \\
-\frac{1}{(m-1)!} & {\bf if}~~\{i_1,i_2,\cdots,i_m\}=V_0\cup V_l \\
0 & {\bf otherwise}.
\end{array}
\right.
$$
It can be verified that $\mathcal{L}{\bf x}^m=\sum_{l=1}^s\mathcal{L}_{\Gamma_l}{\bf x}_{\Gamma_l}^m$ for all ${\bf x}\in \mathbb{R}^n$, and
$\mathcal{L}$ is a symmetric $W$-tensor from Definition \ref{def22}.


\qed

\subsection{Laplacian tensors of hyper-trees}
The notion of hyper-tree is defined based on a normal graph. Let $G=(V, E)$ be a usual graph, that is a 2-uniform hypergraph. If any two
vertices of $G$ are connected by exactly one path, then $G$ is called a tree. In particular, a tree is called a rooted oriented tree if one
vertex has been designated the root, in which case the edges have a natural orientation which one orders the edge from the top of the rooted tree (root)
to the bottom and from the left to the right.

According to \cite{Hu2013}, an $m$-uniform hyper-tree $G=(V, E)$ is the $m$th-power of a tree such that there exists a tree $T=(V_0, E_0)$ and additional vertices
$\bar{V}=\{i_{e,1},i_{e,2},\cdots,i_{e,m-2}~|~ e\in E_0\}$ satisfying $V=V_0 \cup \bar{V}$. If the vertices of $\bar{V}$ are all distinct
and the tree is an oriented rooted tree, we call it an {\bf $m$-uniform oriented rooted hyper-tree generated by independent vertices}. As an illustration, the following hypergraph
$
V=\{1,\ldots,19\}~~{\bf with}~~E=\{(1,2,3,4), (4,5,6,7), (4,8,9,10), (1,11,12,13), (13,14,15,16), (16,17,18,19)\}
$
is a 4-uniform hyper-tree generated by independent vertices, because it can be formed by the tree $T=(V_0, E_0)$, where
$$
V_0=\{1,4,7,10,13,16,19\},~~E_0=\{(1,4),(4,7),(4,10),(1,13),(13,16),(16,19)\}
$$
and the additional vertices $\bar{V}=\{2,3,5,6,8,9,11,12,14,15,17,18\}$.

\bt\label{them42} Let $G=(V, E)$ be an $m$-uniform oriented rooted hyper-tree generated by independent vertices. Then, its Laplacian
tensor is a symmetric $W$-tensor.
\et
\proof Let the adjacent tensor and Laplacian tensor of $G$ be denoted by $\mathcal{A}$ and $\mathcal{L}$ respectively.
Let $\mathcal{D}$ be the diagonal tensor with its diagonal entries $d_i$ i.e. the degree of the vertex $i\in [n]$.
So, it holds that $\mathcal{L}=\mathcal{D}-\mathcal{A}$. From Remark \ref{remark21}, it suffices to show that $-\mathcal{A}$ is
a $W$-tensor. As $G=(V, E)$ is an $m$-uniform oriented rooted hyper-tree by independent vertices, there exists a oriented tree
$T=(V_0, E_0)$ and additional distinct vertices $\bar{V}=\{i_{e,1},i_{e,2},\cdots,i_{e,m-2}~|~ e\in E_0\}$ satisfying $V=V_0 \cup \bar{V}$.
Without loss of generality, suppose
$$
V_0=\{1,2,\cdots,n_0\}~~{\bf and}~~E_0=\{e^1_0,e^2_0,\cdots,e^s_0\},~~{\bf where}~~s=|E_0|.
$$
Let $\bar{\Gamma}_l=\{i\in V_0~|~i\in e^l_0\},~ l\in [s]$. Then, it satisfies that $V_0=\bigcup_{l=1}^s\bar{\Gamma}_l$ and $|\left(\bigcup_{l=1}^{p-1}\bar{\Gamma}_l\right)\bigcap \bar{\Gamma}_p|\leq 1$ for all $2\leq p\leq s$.
Now, define $\Gamma_l= \bar{\Gamma}_l\bigcup \{i_{e_0^l,1},i_{e_0^l,2},\cdots,i_{e_0^l,m-2}\}$, $l\in [s]$. Then, we see that
each $\Gamma_l$ corresponds an edge of the hyper-tree $G$, and it can be easily verified that all conditions {\bf (i)-(iii)} of
Definition \ref{def22} are satisfied. Hence $-\mathcal{A}$ is a $W$-tensor, and it follows that the Laplacian tensor $\mathcal{L}$
is a $W$-tensor. 
\qed

\subsection{Numerical examples}
Throughout this section, all numerical experiments are performed on a desktop, with  3.47 GHz quad-core Intel E5620 Xeon 64-bit CPUs and 4 GB RAM,
equipped with  Matlab 2015.

\begin{example}\label{exam41} Let $G=(V, E)$ be a 4-uniform hyper-star. Suppose $E=\{e_1,e_2,\cdots,e_k\}$, where $k\in \mathbb{N}$.
Then, it holds that $|V|=3k+1$. Without loss of generality, assume
$$
e_j=\{1, 3j-1, 3j, 3j+1\},~j\in [k].
$$
So, the Laplacian tensor $\mathcal{L}$ 
of the hyper-star are 4th order and $(3k+1)$-dimensional tensors such that
$$
\mathcal{L}_{i_1i_2i_3 i_4}=\left\{
\begin{array}{cll}
k & {\bf if}~~i_1=i_2=i_3=i_4=1, \\
1 & {\bf if}~~i_1=i_2=i_3=i_4=i,~i\in\{2,3,\cdots, 3k+1\}, \\
-\frac{1}{3!} & {\bf if}~~(i_1,i_2,i_3,i_4)=(1,3j-1,3j,3j+1),~{\bf for~some}~j\in [k], \\
0 & {\bf otherwise}.
\end{array}
\right.
$$
So, the homogeneous polynomials corresponding to $\mathcal{L}$ is 
$$
\mathcal{L}{\bf x}^4=kx_1^4+x_2^4+\cdots+x_n^4-4\sum_{j=1}^kx_1x_{3j-1}x_{3j}x_{3j+1}.
$$
We now compute the maximum $H$-eigenvalues of $\mathcal{L}$ 
using Algorithm \ref{algo32}.

\begin{table}[!htb]
\begin{center}
\caption{Test results for Laplacian tensor of Example \ref{exam41}}
\medskip
\begin{tabular}{|ccc|c|ccc|}\hline
$m$& $k$  & $n$  & True $\lambda_{\max}^H(\mathcal{L})$ & Est. $\lambda_{\max}^H(\mathcal{L})$ & YALMIP & SeDuMi \\
\hline
 4 & ~~10 & ~~31 & ~~10.0137  & ~~10.0137 & ~~10.1 & ~~3.0  \\
 4 & ~100 & ~301 & ~100.0001  & ~100.0001 & ~~97.4 & ~~5.0  \\
 4 & ~500 & 1501 & ~500.0000  & ~500.0000 & ~557.9 & ~28.6  \\
 4 & 1000 & 3001 & 1000.0000  & 1000.0000 & 1112.1 & ~64.0  \\
 4 & 2000 & 6001 & 2000.0000  & 2000.0000 & 1907.2 & 205.1  \\   \hline
\end{tabular}
\end{center}
\end{table}

Using Algorithm 3.2, we compute the maximum H-eigenvalue of the Laplacian tensor $\lambda_{\max}^H(\mathcal{L})$
of four-uniform hyper-stars with edges ranging from ten to two thousand.
For each case, we obtain $\lambda_{\max}^H(\mathcal{L})$ within $36$ minutes.
It is known that the true maximum $H$-eigenvalue $\lambda_{\max}^H(\mathcal{L})$ is the unique root of the polynomial equation $(1-x)^{m-1}(x-k)+k=0$ in the open interval $(k,k+1)$, where $k$ is the number
of edges and $m$ is the degree of the hypergraph.
The preceding polynomial equation is solved by the Matlab command ``${\rm vpnsolve}$''.
The results are summarized in Table 2 where the meanings of the data are the same as in Table 1. It can be easily seen that the maximum H-eigenvalue estimated by Algorithm 3.2 are consistent with the true maximum H-eigenvalues.

%
%
%
\end{example}


\begin{example}\label{example63} Let $k, m\in \mathbb{N}$. Suppose $m$ is an even number. Assume $G=(V, E)$ is an $m$-uniform hypergraph with vertices and edges such that
$$
V=\{1,2,\cdots,k(m-1)+1\},~~E=\{e_1,e_2,\cdots,e_k\},
$$
where each edge $e_l=\{(l-1)(m-1)+1,(l-1)(m-1)+2,\cdots, l(m-1)+1\}$, $l\in [k]$.
So, the hypergraph $G=(V,E)$ is the case of a hyper-tree with its Laplacian tensor $\mathcal{L}$ is an $m$th order $n=k(m-1)+1$ dimensional tensor.
By a direct computation, we obtain
$$
\mathcal{L}_{i_1i_2\cdots i_m}=\left\{
\begin{array}{cll}
2 & {\bf if}~~i_1=i_2=\cdots=i_m=l(m-1)+1,~l\in \{1,2,\cdots,k-1\}, \\
1 & {\bf if}~~i_1=i_2=\cdots=i_m=i,~i\in [n]\backslash \{(m-1)+1, 2(m-1)+1,\cdots,(k-1)(m-1)+1\},  \\
-\frac{1}{(m-1)!} & {\bf if}~~\{i_1,i_2,\cdots,i_m\}=e_l,~{\bf for~some}~l\in [k], \\
0 & {\bf otherwise}.
\end{array}
\right.
$$
Moreover, we have that
$$
\mathcal{L}{\bf x}^m=\sum_{i\in [n]\backslash \{l(m-1)+1|l\in [k-1]\}} x_i^m+2\sum_{l=1}^{k-1}x^m_{l(m-1)+1}-m\sum_{l=1}^kx_{(l-1)(m-1)+1}x_{(l-1)(m-1)+2}\cdots x_{l(m-1)+1}.
$$
Combining this with Algorithm \ref{algo32}, we compute the maximum $H$-eigenvalues of $\mathcal{L}$ 
with different $m$ and $n$, and the results are listed in the Table 3.
It should be noted that Hu et al. \cite{Hu2015} showed that the maximum H-eigenvalues of the Laplacian tensor $\mathcal{L}$
and the signless Laplacian tensor $\mathcal{Q}$ are equivalent when the even-uniform hypergraph is connected and odd-bipartite.
On the other hand, since the signless Laplacian tensor is nonnegative, its maximum $H$-eigenvalue could be computed
by the classical NQZ algorithm \cite{NQZ2009}. Direct verification shows that the hypergraph discussed in this example is connected and odd-bipartite, and hence, the maximum $H$-eigenvalue of its Laplacian tensor $\mathcal{L}$ can also be computed by the NQZ method. To compare the performance of our method and the
NQZ method, we list the estimated maximum H-eigenvalue of the signless Laplacian tensor $\lambda_{\max}^H(\mathcal{Q})$
as well as CPU-time of the NQZ algorithm in Table 3. It can be seen that estimated maximum H-eigenvalues of $\lambda_{\max}^H(\mathcal{L})$ and
$\lambda_{\max}^H(\mathcal{Q})$ coincide which numerically verifies the assertion in \cite{Hu2015}.
Interestingly, one also observes that Algorithm 3.2 is indeed {\em faster than} the first-order NQZ algorithm \cite{NQZ2009} as it exploits
the structure of the underlying problem.
%

\begin{table}[h!]
\begin{center}
\caption{Test results for Laplacian tensor of Example \ref{example63}}
\medskip
\begin{tabular}{|ccc|ccc|cc|}\hline
$m$ & $k$  & $n$& Est. $\lambda_{\max}^H(\mathcal{L})$ & YALMIP & SeDuMi  & Est. $\lambda_{\max}^H(\mathcal{Q})$ & NQZ \\   \hline
 4  & ~100 & ~301 & 2.9997 & ~~~4.4 & ~~3.3 & 2.9997 & ~~~9.4 \\
 4  & ~400 & 1201 & 3.0000 & ~~28.4 & ~12.7 & 3.0000 & ~337.0 \\
 4  & 1000 & 3001 & 3.0000 & ~160.6 & ~36.4 & 3.0000 & 2678.3 \\   \hline
 6  & ~100 & ~501 & 2.6954 & ~~19.9 & ~14.2 & 2.6954 & ~~24.3 \\
 6  & ~400 & 2001 & 2.6956 & ~282.1 & ~61.7 & 2.6956 & ~923.9 \\
 6  & 1000 & 5001 & 2.6956 & 2828.3 & 144.2 & 2.6956 & 6959.4 \\   \hline
\end{tabular}
\end{center}
\end{table}
\end{example}

\section{Applications in copositivity test of tensors}

In this section, we present further applications on testing the copositivity of a multivariate form associated with symmetric extended
Z-tensors, where the order of the tensor can be either odd or even.  We first recall the definition of extended $Z$-tensors \cite{chen15}.
\bd\label{def61}
A symmetric tensor $\mathcal{A}$ is called an {\bf extended $Z$-tensor} if its associated polynomial $f_{\mathcal{A}}({\bf x})=\mathcal{A}{\bf x}^m$
satisfies that there exist $s \in \mathbb{N}$ with $s \le n$ and index sets $\Gamma_l \subseteq \{1,\cdots,n\}$, $l=1,\cdots,s$ with  $\bigcup_{l=1}^s \Gamma_l=\{1,\cdots,n\}$ and $\Gamma_{l_1} \cap \Gamma_{l_2} =\emptyset$ for all $l_1 \neq l_2$  such that
\[
f({\bf x})=\sum_{i=1}^n f_{m,i} x_i^{m}+\sum_{l=1}^s \sum_{\alpha_l \in \Omega_l }f_{\alpha_l} {\bf x}^{\alpha_l},\]
where
\[
\Omega_l=\left\{\alpha\in ([n]\cup \{0\})^n:~|\alpha|=m, {\bf x}^{\alpha}=x_{i_1}x_{i_2}\cdots x_{i_m}, \{i_1,\cdots,i_m\}\subseteq \Gamma_l,~
\mbox{and} ~~ \alpha \neq m {\bf e_i}, \ i=1,\cdots,n \right \}
\]
for each $l=1,\cdots,s$ and either one of the following two conditions holds:
\begin{itemize}
\item[{\rm (1)}] $f_{\alpha_l}=0$ for all but one $\alpha_l \in \Omega_l$;
\item[{\rm (2)}] $f_{\alpha_l} \le 0$ for all $\alpha_l \in \Omega_l$.
\end{itemize}
\ed

Let $\mathcal{A}=(a_{i_1i_2\cdots i_m})$ be a symmetric tensor with order $m$ and dimension $n$.
Then $\mathcal{A}$ is copositive if and only if
$$
\mathcal{A}{\bf x}^m=\sum_{i_1,\ldots,i_m=1}^n a_{i_1 \ldots i_m}x_{i_1} \cdots x_{i_m}\geq0,~\forall~{\bf x}\in \mathbb{R}^n_+,
$$
which is equivalent to
\begin{equation}\label{e62}
h({\bf x})=\mathcal{A}_h{\bf x}^{2m}=\sum_{i_1,\ldots,i_m=1}^n a_{i_1 \ldots i_m}x_{i_1}^2 \cdots x_{i_m}^2\geq0,~\forall~{\bf x}\in \mathbb{R}^n,
\end{equation}
where $\mathcal{A}_h$ is a symmetric tensor with order $2m$ and dimension $n$.

In particular, if $\mathcal{A}$ is a symmetric extended $Z$-tensor (odd or even order), then $\mathcal{A}_h$ is also an even order extended $Z$-tensor \cite{chen15}. Thus, $-\mathcal{A}_h$ is an even order $W$-tensor. Let $f({\bf x})=-\mathcal{A}_h{\bf x}^{2m}$.
Then, by Theorem \ref{them32} and (\ref{e62}), we have the following corollary.
\bc\label{corol61}
Let $\mathcal{A}$ be a symmetric extended $Z$-tensor with order $m$ and dimension $n$.
For ${\bf x}\in \mathbb{R}^n$, suppose $\mathcal{A}_h$ and $f({\bf x})$ are defined as above.
Then, $\mathcal{A}$ is copositive if and only if
\begin{equation}\label{e63}
\min_{t\in \mathbb{R}} \{t~|~t\|{\bf x}\|_{2m}^{2m}-f({\bf x})\in \Sigma_{2m}^2[{\bf x}]\}\leq 0.
\end{equation}
\ec

We now use the above corollary to test the copositivity of symmetric extended $Z$-tensors with order $m$ and dimension $n$.
The concrete process is listed below.
\bigskip

\newpage

{\bf Procedure }
\begin{itemize}
\item[(i)] Given $(m, n,s,k,M)$ with $n=sk$, where $n$ and $m$ are the dimension and the order of the
randomly generated tensor, respectively,  and $M$ is a large positive constant.

\item[(ii)] Randomly generate a partition of the index set $\{1,\cdots,n\}$, $\{\Gamma_1,\cdots,\Gamma_s\}$, such that $|\Gamma_i|=k$, $i=1,\cdots,s$ and $\Gamma_{i} \cap \Gamma_{i'}=\emptyset$ for all $i \neq i'$.
For each $i=1,\cdots,s-1$, generate a random multi-index $(l_1^i,\cdots,l_m^i)$ with $l_j^i \in \Gamma_i$, $j=1,\cdots,m$ and a random number $\bar{a}_{l_1^i \cdots l_m^i} \in [0,1]$.  { Generate one randomly $m$-th-order $k$-dimensional symmetric tensor $\mathcal{B}$,  such that all elements
of $\mathcal{B}$ are in the interval $[0, 1]$.}

\item[(iii)]  We define extended $Z$-tensor $\mathcal{A}=(a_{i_1i_2\cdots i_m})$ such that
\[
a_{i_1 \cdots i_m}= \left\{\begin{array}{cll}
 M  & \mbox{ if } & i_1=\cdots=i_m=i \mbox{ for all } i=1,\cdots,n, \\
\bar{a}_{l_1^i \cdots l_m^i} & \mbox{ if } & (i_1,\cdots,i_m)=\sigma(l_1^i,\cdots,l_m^i) \mbox{ with } l_1^i,\cdots,l_m^i \in \Gamma_i, i\in [s-1], \\
-\mathcal{B}_{i_1 \cdots i_m} & \mbox{ if } & {  i_1,\cdots,i_m \in \Gamma_s,}\\
0 & \mbox{ otherwise. } &
\end{array}
\right.
\]
Here $\sigma(i_1,\cdots,i_m)$ denotes all the possible permutation of $(i_1,\cdots,i_m)$.

\item[(iv)] Let $\mathcal{A}_h=(a^h_{i_1i_2\cdots i_{2m}})$ be a extended $Z$-tensor with order $2m$ and dimension $n$ such that
$$a^h_{\sigma(i_1i_1i_2i_2\cdots i_mi_m)}=a_{i_1i_2\cdots i_m},~\forall~i_1,i_2,\cdots,i_m\in [n],$$
and $a^h_{i_1i_2\cdots i_{2m}}=0$ otherwise.

\item[(v)] Let $f({\bf x})=-\mathcal{A}_h{\bf x}^{2m}$, ${\bf x}\in \mathbb{R}^n$. Then solve the
SOS programming problem (\ref{e63}) by Matlab Toolbox {\rm YALMIP} \cite{Yalmip1,Yalmip2}
and SeDuMi \cite{Sturm99}.
\end{itemize}

We perform the above procedure to test extended Z-tensors with order $m=3,4,5,6$
and dimension $n=2500,2000,1500,1000$, respectively.
For each case, we fix $s=500$ and run one hundred tests.
Table \ref{Tab-6.2} summarizes the percentage of copositive instances of these extended $Z$-tensors.
Clearly, for fixed order $m$ and dimension $n$, the percentage of copositive extended $Z$-tensors increase
as the parameter $M$ increase.
Moreover, if $M$ is large enough, the generated extended $Z$-tensor must be positive definite, and so, in particular is copositive.

\begin{table}[bth!]
  \caption{The percentage of copositive instances of randomly generated extended $Z$-tensors.}\label{Tab-6.2}
  \centering
  \begin{tabular}{|c|cccccc|}
    \hline
    \multicolumn{7}{|c|}{$m=3$ and $n=2500$} \\ \hline
    $M$          &  11  &  12  &  13  &  14  &  15  &  16  \\ \hline
    Copositivity & 12\% & 35\% & 68\% & 91\% & 96\% & 99\% \\
    \hline\hline
    \multicolumn{7}{|c|}{$m=4$ and $n=2000$} \\ \hline
    $M$          &  25  &  28  &  31  &  34  &  37  &  40  \\ \hline
    Copositivity &  2\% & 11\% & 36\% & 63\% & 89\% & 99\% \\
    \hline\hline
    \multicolumn{7}{|c|}{$m=5$ and $n=1500$} \\ \hline
    $M$          &  30  &  35  &  40  &  45  &  50  &  55  \\ \hline
    Copositivity &  6\% & 19\% & 50\% & 74\% & 93\% & 99\% \\
    \hline\hline
    \multicolumn{7}{|c|}{$m=6$ and $n=1000$} \\ \hline
    $M$          &   9  &  12  &  15  &  18  &  21  &  24  \\ \hline
    Copositivity &  5\% & 20\% & 44\% & 67\% & 87\% & 99\% \\
    \hline
  \end{tabular}
\end{table}

\setcounter{equation}{0}
\section{Conclusions and remarks}
In this article, we propose an efficient semi-definite program algorithm to compute the maximum $H$-eigenvalues of
even order symmetric $W$-tensors  based on its potential SOS structure.
Furthermore, we present two interesting applications:  examining the copositivity of symmetric extended $Z$-tensors with even or odd orders,
and computing the maximum $H$-eigenvalues of Laplacian tensors of hyper-stars and hyper-trees.
Our results point out some interesting future directions. For example, note that  the $W$-tensors considered in this paper are with even order. Can
we compute the maximum $H$-eigenvalues of odd order symmetric $W$-tensors? Can we compute the minimum $H$-eigenvalues of $W$-tensors?
These will be considered in a future work.

\end{document}